\documentclass[preprint,12pt]{elsarticle}

\usepackage{amssymb}

\newtheorem{theorem}{Theorem}[section]
\newtheorem{lemma}[theorem]{Lemma}
\newtheorem{proposition}[theorem]{Proposition}
\newtheorem{corollary}[theorem]{Corollary}

\newtheorem{question}[theorem]{Question}
\newtheorem{definition}[theorem]{Definition}
\newtheorem{example}[theorem]{Example}

\journal{}

\begin{document}

\begin{frontmatter}



\title{On surjectively universal Polish groups}


\author{Longyun Ding\footnote{Research partially supported by the National Natural Science Foundation of China (Grant No. 10701044)
and Program for New Century Excellent Talents in University.}}

\ead{dinglongyun@gmail.com}
\address{School of Mathematical Sciences and LPMC, Nankai University, Tianjin 300071, PR China}

\begin{abstract}

A Polish group is surjectively universal if it can be continuously homomorphically mapped onto every Polish group.
Making use of a type of new metrics on free groups \cite{DG}, we prove the existence of surjectively universal Polish groups,
answering in the positive a question of Kechris. In fact, we give several examples of surjectively universal Polish groups.

We find a sufficient condition to guarantee that the new metrics on free groups can be computed directly. We also compare this condition with CLI groups.
\end{abstract}


\end{frontmatter}



\section{Introduction}

A topological group is called a {\it Polish group} if its underlying space is a Polish space, i.e. a separable, completely metrizable topological space.
For Polish groups, there are usually two parallel notions of universality.
A Polish group is {\it universal} if every Polish group is isomorphic to one of its closed subgroup. We call a Polish group $G$ {\it surjectively universal}
if every Polish group $H$ is isomorphic to a topological quotient group of $G$, or equivalently, there is a continuous surjective homomorphism $\Phi:G\to H$.

It was proved by Uspenski\u{\i} \cite{uspenskii} that there exists a universal Polish group.
For surjectively universal, Kechris asked the following question (Problem 2.10 of \cite{kechris}).
\begin{question}\label{main}
Is there a surjectively universal Polish group?
\end{question}
This question arose again in Becker-Kechris' book (Open problem 1.4.2 of \cite{BK}), and was also mentioned in \cite{gaobook}.

Some related problems were investigated in these years. Some of them are on universality of subclasses of Polish groups.
Recall that a metric $d$ on a group $G$ is {\it left-invariant} if $d(gh,gk)=d(h,k)$ for all $g,h,k\in G$. The definition of right-invariant metric is similar.
$d$ is {\it two-sided invariant} if it is both left- and right-invariant. A Polish group is {\it CLI} if it admits a compatible complete left-invariant metric.
Endowed with the Graev metrics \cite{graev} on free groups, it was shown in \cite{SPW} that there is a surjectively universal in the class of all Polish groups
which admit an two-sided invariant compatible complete metric. This implies the existence of a surjectively universal abelian Polish group.
It is also known that there are universal Polish abelian groups (see \cite{shkarin}). The notion of weakly universal is introduced in \cite{GP}.
A Polish group $G$ is {\it weakly universal} if every Polish group is isomorphic to a closed subgroup of a topological quotient group of $G$.
For example, $\ell_1$ under addition is a weakly universal for abelian Polish groups. On the negative side,
most recently, Malicki \cite{malicki} proved that there is no weakly universal CLI group.
So neither universal CLI group, nor surjectively universal CLI group exists.

With a notion of {\it scale}, a type of new metrics on free groups were constructed in \cite{DG}.
As a result, a class of Polish groups were obtained by taking completion of free groups with the new metrics such that,
every Polish group is isomorphic to a topological quotient group of one in the class. Thus if there exist surjectively universal Polish groups,
there shall be one in this class. Furthermore, it implies that,
every Polish group is isomorphic to a topological quotient group of a ${\bf\Pi}^0_3$ subgroup of $S_\infty$.
In \cite{tsankov}, another class of ${\bf\Pi}^0_3$ subgroups of $S_\infty$ with similar surjective universality properties were obtained by a very different way.

In this article, we give a sufficient condition on scales to guarantee that the completion of free groups with the induced metrics are surjectively universal.
Granted with the condition, we give several examples of surjectively universal Polish groups, answering in the positive of Question~\ref{main}.

In general, the value of metrics on free groups defined in \cite{DG} are very hard to compute. We also find a suitable subclass of scales,
which will be named {\it adequate} scales, such that the induced metrics can be computed directly. Another example of surjectively universal Polish group,
in which the involved scale is adequate, is also given. Furthermore, this example can be isomorphic to a ${\bf\Pi}^0_3$ subgroup of $S_\infty$.
In the end, we will show that any non-trivial adequate scale can not induce a CLI group.

\section{Review of metrics on free groups}

In this section we review definitions and several results on metrics on free groups defined in \cite{DG}.
These definitions and results will play a central role throughout the rest of this paper.

For a nonempty set $X$, we define the free group on $X$ slightly different from the usual definition. Let $X^{-1}=\{x^{-1}:x\in X\}$ be a disjoint copy of $X$,
and let $e\notin X\cup X^{-1}$. We will use $e$ rather than the empty word to denote the identity element of free group.
We use notation convention that $(x^{-1})^{-1}=x$ for $x\in X$ and $e^{-1}=e$.

Put $\overline{X}=X\cup X^{-1}\cup\{e\}$. Let $W(X)$ be the set of words on $\overline{X}$. For $w\in W(X)$, ${\rm lh}(w)$ stands for the length of $w$.
A word $w\in W(X)$ is {\it irreducible} if $w=e$ or $w=x_0\cdots x_n$ with $x_i\ne e$ and $x_{i+1}\ne x_i^{-1}$ for each $i$.
Let $F(X)$ be the set of all irreducible words. For each $w\in W(X)$,
the {\it reduced word} for $w$, denoted $w'$, is the unique irreducible word obtained by successively replacing any occurrence of $xx^{-1}$ in $w$ by $e$,
and eliminating $e$ from any occurrence of the form $w_1ew_2$, where at least one of $w_1$ and $w_2$ is nonempty.
We also say $w$ is a {\it trivial extension} of $w'$. By defining $w\cdot u=(wu)'$,
where $wu$ is the concatenation of $w$ and $u$, we turn $F(X)$ into a group, i.e. the {\it free group} on $X$.

Assume now $(X,d)$ is a metric space. We extend $d$ to be a metric on $\overline{X}$, still denoted by $d$, such that, for all $x,y\in\overline{X}$, we have
$d(x,y)=d(x^{-1},y^{-1})$. For further extending $d$ to be a metric on the free group $F(X)$, we need two more notation.

\begin{definition}[Ding-Gao, \cite{DG}]{\rm
Let $\Bbb R_+$ denote the set of non-negative real numbers. A function $\Gamma:\overline{X}\times\Bbb R_+\to\Bbb R_+$ is a {\it scale} on $\overline{X}$ if
the following hold for any $x\in\overline{X}$ and $r\in\Bbb R_+$:
\begin{enumerate}
\item[(i)] $\Gamma(e,r)=r$; $\Gamma(x,r)\ge r$;
\item[(ii)] $\Gamma(x,r)=0$ iff $r=0$;
\item[(iii)] $\Gamma(x,\cdot)$ is a monotone increasing function with respect to the second variable;
\item[(iv)] $\lim_{r\to 0}\Gamma(x,r)=0$.
\end{enumerate}}
\end{definition}

Let $G$ be a metrizable group and $d_G$ be a compatible left-invariant metric on $G$. For $g\in G$ and $r\ge 0$, we denote
$$\Gamma_G(g,r)=\max\{r,\sup\{d_G(1_G,g^{-1}hg):d_G(1_G,h)\le r\}\},$$
where $1_G$ is the identity element of $G$.

\begin{proposition}\label{gammag}
For $g_1,g_2\in G$ and $r\ge 0$, we have
$$|\Gamma_G(g_1,r)-\Gamma_G(g_2,r)|\le 2d_G(g_1,g_2).$$
\end{proposition}

{\bf Proof.} Note that $d_G$ is left-invariant. For $h\in G$, if $d(1_G,h)\le r$, then
$$\begin{array}{ll}d_G(1_G,g_1^{-1}hg_1)&\le d_G(1_G,g_1^{-1}g_2)+d_G(g_1^{-1}g_2,g_1^{-1}hg_2)+d_G(g_1^{-1}hg_2,g_1^{-1}hg_1)\cr
&=d_G(g_1,g_2)+d_G(1_G,g_2^{-1}hg_2)+d_G(g_2,g_1).\end{array}$$
It follows that $\Gamma_G(g_1,r)\le\Gamma_G(g_2,r)+2d_G(g_1,g_2)$.
\hfill$\Box$

\begin{definition} [Ding-Gao, \cite{DG1}]{\rm
Let $m,n\in\Bbb N$ and $m\le n$. A bijection $\theta$ on $\{m,\cdots,n\}$ is a {\it match} if
\begin{enumerate}
\item[(1)] $\theta\circ\theta={\rm id}$; and
\item[(2)] there is no $m\le i,j\le n$ such that $i<j<\theta(i)<\theta(j)$.
\end{enumerate}}
\end{definition}

Before defining a norm and a metric on $F(X)$ with a given scale $\Gamma$, we must define a pre-norm on $W(X)$ firstly.

\begin{definition}[Ding-Gao, \cite{DG}]\label{prenorm}{\rm
Let $\Gamma$ be a scale on $\overline{X}$. For $l\in\Bbb N$, $w\in W(X)$ with ${\rm lh}(w)=l+1$ and $\theta$ a match on $\{0,\cdots,l\}$,
we define {\it pre-norm} $N_\Gamma^\theta(w)$ by induction on $l$ as follows:
\begin{enumerate}
\item[(0)] if $l=0$, let $w=x$ and define $N_\Gamma^\theta(w)=d(e,x)$; else if $l=1$ and $\theta(0)=1$, let $w=x^{-1}y$ and define $N_\Gamma^\theta(w)=d(x,y)$;
\item[(1)] if $l>0$ and $\theta(0)=k<l$, let $\theta_1=\theta\upharpoonright\{0,\cdots,k\}$, $\theta_2=\theta\upharpoonright\{k+1,\cdots,l\}$
and $w=w_1w_2$ where ${\rm lh}(w_1)=k+1$; define
$$N_\Gamma^\theta(w)=N_\Gamma^{\theta_1}(w_1)+N_\Gamma^{\theta_2}(w_2);$$
\item[(2)] if $l>1$ and $\theta(0)=l$, let $\theta_1=\theta\upharpoonright\{1,\cdots,l-1\}$ and $w=x^{-1}w_1y$ where $x,y\in\overline{X}$;
then ${\rm lh}(w_1)=l-1$ and $w=(x^{-1}w_1x)\cdot(x^{-1}y)=(x^{-1}y)\cdot(y^{-1}w_1y)$. Define
$$N_\Gamma^\theta(w)=d(x,y)+\min\{\Gamma(x,N_\Gamma^{\theta_1}(w_1)),\Gamma(y,N_\Gamma^{\theta_1}(w_1))\}.$$
\end{enumerate}}
\end{definition}

\begin{definition}[Ding-Gao, \cite{DG}]{\rm
Given a scale $\Gamma$ on $\overline{X}$, we define a {\it norm} $N_\Gamma$ and a metric $\delta_\Gamma$ on $F(X)$ by
$$N_\Gamma(w)=\inf\{N_\Gamma^\theta(w^*):(w^*)'=w,\theta\mbox{ is a match}\},$$
$$\delta_\Gamma(w,v)=N_\Gamma(w^{-1}\cdot v).$$}
\end{definition}

The metric $\delta_\Gamma$ is left-invariant. If $X$ is separable, then $F(X)$ is a separable topological group with respect to
the topology induced by $\delta_\Gamma$ (see \cite{DG}, Theorem 3.9). In the next section, we will use the following lemma, i.e. Lemma 3.7 of \cite{DG},
to prove the existence of surjectively universal Polish groups.

\begin{lemma}[Ding-Gao]\label{extend}
Let $G$ be a topological group and $d_G$ a compatible left-invariant metric on $G$. Let $\Gamma$ be a scale on $\overline{X}$.
Let $\varphi:\overline{X}\to G$ be a function. Suppose that for any $x,y\in\overline{X}$ and $r\ge 0$:
\begin{enumerate}
\item[(i)] $\varphi(e)=1_G$; $\varphi(x^{-1})=\varphi(x)^{-1}$;
\item[(ii)] $d_G(\varphi(x),\varphi(y))\le d(x,y)$; and
\item[(iii)] $\Gamma_G(\varphi(x),r)\le\Gamma(x,r)$.
\end{enumerate}
Then $\varphi$ can be uniquely extended to a group homomorphism $\Phi:F(X)\to G$ such that, for any $w,v\in F(X)$,
$$d_G(\Phi(w),\Phi(v))\le\delta_\Gamma(w,v).$$
\end{lemma}

For $w,u\in F(X)$, let
$$\Delta_\Gamma(w,u)=\max\{\delta_\Gamma(w,u),\delta_\Gamma(w^{-1},u^{-1})\}.$$
Then $\Delta_\Gamma$ is a compatible metric on $F(X)$. Denote by $\overline{F}_\Gamma(X)$ the completion of $(F(X),\Delta_\Gamma)$.
Then $\overline{F}_\Gamma(X)$ is a Polish group. Furthermore, $\delta_\Gamma$ can be uniquely extend onto $\overline{F}_\Gamma(X)$ which is still a compatible
left-invariant metric.

\section{Surjectively universal Polish groups}

The objective of this section is to find a sufficient condition on the scale $\Gamma$ such that $\overline{F}_\Gamma(X)$ is surjectively universal.
Here the original metric space $X$ shall be a special subset of the Baire space ${\cal N}=\Bbb N^\Bbb N$. For distinct $x,y\in{\cal N}$, we define
$$d(x,y)=\max\{2^{-n}:x(n)\ne y(n)\},$$
$$d(x^{-1},y^{-1})=d(x,y),\quad d(x,e)=d(x^{-1},e)=d(x,y^{-1})=1.$$
Then $d$ is metric on $\overline{\cal N}$. For $n\in\Bbb N$, let
$${\cal N}_n=\{x\in{\cal N}:\forall m\ge n(x(m)=0)\}.$$
$$\pi_n(x)(m)=\left\{\begin{array}{ll}x(m), &\mbox{ if }m<n,\cr 0, &\mbox{ if }m\ge n.\end{array}\right.$$
And let ${\cal N}_\omega=\bigcup_{n\in\Bbb N}{\cal N}_n$. It is easy to see that ${\cal N}_\omega$ is dense in $\cal N$.

Denote $\hat 0=\langle 0,0,\cdots\rangle$ the unique element in ${\cal N}_0$.

We now focus on scales on $\overline{\cal N}_\omega$.

\begin{theorem}\label{universal}
Let $\Gamma$ be a scale on $\overline{\cal N}_\omega$ such that $\Gamma(x,r)=\Gamma(x^{-1},r)$ for any $x\in{\cal N}$.
If there exists a $K\in\Bbb N$ satisfying that, for $m>0$ and $k>m+K+5$,
\begin{enumerate}
\item[(S1)] for any $x\in{\cal N}_m\setminus{\cal N}_{m-1}$ and $r>2^{-(m+K+5)}$,
there exists $y\in{\cal N}_{m-1}$ such that $\Gamma(x,r)\ge\min\{2^{-(K+3)},\Gamma(y,r)+2^{-K}d(x,y)\}$;
\item[(S2)] by defining subset $E^k_m\subseteq{\cal N}_m\setminus{\cal N}_{m-1}$ as
$$\begin{array}{ll}x\in E^k_m\iff &\exists r\in(2^{-k},2^{-(m+K+5)}]\forall y\in{\cal N}_{m-1}\cr
&(\Gamma(x,r)<\min\{2^{-(K+3)},\Gamma(y,r)+2^{-K}d(x,y)\}),\end{array}$$
then $E^k_m$ is finite; and
\item[(S3)] $\Gamma(x,r)\ge 8r$ for any $x\in{\cal N}_m\setminus{\cal N}_{m-1}$ and $r\le 2^{-(m+K+5)}$.
\end{enumerate}
Then $\overline{F}_\Gamma(\cal N_\omega)$ is a surjectively universal Polish group.
\end{theorem}

{\bf Proof.}
Let $G$ be a Polish group, and let $(U_n)_{n\in\Bbb N}$ be a neighborhood base of its identity element $1_G$ such that $U_0=G$, and
$U_n=U_n^{-1}$, $U_{n+1}^3\subseteq U_n$ for all $n\in\Bbb N$ (see, e.g. \cite{gaobook}, Theorem 2.2.1).

Firstly, we inductively define a strictly increasing $f:\Bbb N\to\Bbb N$ and a $\varphi:{\cal N}_\omega\to G$.
We put $f(0)=0,f(1)=1$ and $\varphi(\hat 0)=1_G$.
Let $n\ge 1$. Suppose we have defined $f(n)$ and $\varphi:{\cal N}_{n-1}\to G$.
We shall extend $\varphi$ onto ${\cal N}_n$ such that
\begin{enumerate}
\item[(a)] $\forall x\in{\cal N}_n(\varphi(x)\in[\varphi(\pi_{n-1}(x))U_{f(n-1)}\cap U_{f(n-1)}\varphi(\pi_{n-1}(x))])$; and
\item[(b)] for each $y\in{\cal N}_{n-1}$ we have
$$[\varphi(y)U_{f(n-1)}\cap U_{f(n-1)}\varphi(y)]\subseteq
\bigcup_{\scriptsize\begin{array}{cc}\pi_{n-1}(x)=y\cr x\in{\cal N}_n\end{array}}[\varphi(x)U_{f(n)}\cap U_{f(n)}\varphi(x)].$$
\end{enumerate}
In order to realize such an extension, we fix a $y\in{\cal N}_{n-1}$. Let $O$ be the open set $\varphi(y)U_{f(n-1)}\cap U_{f(n-1)}\varphi(y)$,
the following is an open cover of $O$:
$$\{gU_{f(n)}\cap U_{f(n)}g:g\in O\}.$$
Since the topology of $G$ is second countable, we can find a countable subcover. Thus there is a countable subset $C\subseteq O$ such that
$$O\subseteq\bigcup_{g\in C}[gU_{f(n)}\cap U_{f(n)}g].$$
Then we extend $\varphi$ such that $\{\varphi(x):x\in{\cal N}_n,\pi_{n-1}(x)=y\}=C$.

To define $f(n+1)$, denote
$$B_n=\bigcap_{m=1}^{n}\bigcap_{x\in E^{n+K+6}_m}[\varphi(x)U_{f(n)}\varphi(x)^{-1}\cap\varphi(x)^{-1}U_{f(n)}\varphi(x)].$$
Since every $E_m^{n+K+6}$ is finite, $B_n$ is an open neighborhood of $1_G$. We can find a sufficiently large $N$ such that $U_N\subseteq B_n$.
Define $f(n+1)=N$. This finishes the induction stage.

We extend $\varphi:{\cal N}_\omega\to G$ onto ${\overline{\cal N}}_\omega$ as
$$\varphi(e)=1_G,\quad\varphi(x^{-1})=\varphi(x)^{-1}.$$

Secondly, we define a compatible left-invariant metric on $G$.

Denote $V_k=U_{f(k)}$ for $k\in\Bbb N$. For $g,h\in G$, let
$$\rho(g,h)=\min\{2^{-(k+K+3)}:g^{-1}h\in V_k\}.$$
By a standard method, we define metric $d_G$ on $G$ as
$$d_G(g,h)=\inf\left\{\sum_{i=0}^l\rho(g_i,g_{i+1}):g_0=g,g_{l+1}=h,g_1,\cdots,g_l\in G,l\in\Bbb N\right\}.$$
Following the proof of Birkhoff-Kakutani's theorem (see, e.g. \cite{gaobook}, Theorem 2.2.1), for $g,h\in G$, we have
$$\frac{1}{2}\rho(g,h)\le d_G(g,h)\le\rho(g,h)\le 2^{-(K+3)}.$$
Moreover, $d_G$ is a compatible left-invariant metric on $G$.

Now we are going to apply Lemma~\ref{extend} on these $\Gamma$ and $d_G$. We need to check clauses (ii) and (iii) of Lemma~\ref{extend}.

For (ii), we let $x,y\in{\cal N}_\omega$ with $d(x,y)=2^{-n}$, then $\pi_n(x)=\pi_n(y)$. Note that
$$d_G(\varphi(x),\varphi(\pi_n(x)))\le\sum_{m\ge n}d_G(\varphi(\pi_{m+1}(x)),\varphi(\pi_m(x))).$$
By (a), we have $\varphi(\pi_{m+1}(x))\in\varphi(\pi_m(x))U_{f(m)}$, so
$$\varphi(\pi_m(x))^{-1}\varphi(\pi_{m+1}(x))\in U_{f(m)}=V_m.$$
Thus $d_G(\varphi(\pi_{m+1}(x)),\varphi(\pi_m(x)))\le 2^{-(m+K+3)}$. It follows that
$$\begin{array}{ll}d_G(\varphi(x),\varphi(y))&\le d_G(\varphi(x),\varphi(\pi_n(x)))+d_G(\varphi(y),\varphi(\pi_n(y)))\cr
&\le 2\sum_{m\ge n}2^{-(m+K+3)}\cr
&<2^{-(n+K+1)}=2^{-(K+1)}d(x,y).\end{array}$$
By the same arguments, we have $d_G(\varphi(x^{-1}),\varphi(y^{-1}))\le 2^{-(K+1)}d(x^{-1},y^{-1})$.

For (iii), we inductively prove that, for $r>0$ and $x\in{\cal N}_m\setminus{\cal N}_{m-1}$, we have $\Gamma_G(\varphi(x),r)\le\Gamma(x,r)$.
If $x=\hat 0$, the inequality is trivial. Note that $\Gamma_G(\varphi(x),r)\le 2^{-(K+3)}$.
We see the desired inequality holds whenever $\Gamma(x,r)\ge 2^{-(K+3)}$. Without loss of generality, we may assume that $\Gamma(x,r)<2^{-(K+3)}$.
Let $2^{-k}<r\le 2^{-(k-1)}$ for some $k\in\Bbb Z$.

If $k\le m+K+5$ or $x\notin E^k_m$, then from (S1) and (S2), there exists $y\in{\cal N}_{m-1}$ such that $\Gamma(x,r)\ge\Gamma(y,r)+2^{-K}d(x,y)$.
By induction hypothesis and Proposition~\ref{gammag}, we have
$$\Gamma(x,r)\ge\Gamma_G(\varphi(y),r)+2d_G(\varphi(x),\varphi(y))\ge\Gamma_G(\varphi(x),r).$$

Otherwise, we have $k>m+K+5$ and $x\in E^k_m$. For any $h\in G$ with $d_G(1_G,h)\le r$, we have $\rho(1_G,h)\le 2r\le 2^{-(k-2)}$, so $h\in V_{k-K-5}$.
By the definition of $f(k-K-5)$, we have $V_{k-K-5}\subseteq B_n$ with $n=k-K-6$, so $m\le n$. Then $\varphi(x)^{-1}h\varphi(x)\in V_{k-K-6}$.
It follows from (S3) that
$$\sup\{d_G(1_G,\varphi(x)^{-1}h\varphi(x)):d(1_G,h)\le r\}\le 2^{-(k-3)}<8r\le\Gamma(x,r).$$
Thus $\Gamma_G(\varphi(x),r)\le\Gamma(x,r)$.

By the same arguments, we can prove that $\Gamma_G(\varphi(x^{-1}),r)\le\Gamma(x^{-1},r)$.

Therefore, Lemma~\ref{extend} applies in this case. $\varphi$ can be extended to a group homomorphism $\Phi:F({\cal N}_\omega)\to G$ such that for any
$w,v\in F({\cal N}_\omega)$
$$d_G(\Phi(w),\Phi(v))\le\delta_\Gamma(w,v).$$
Let $D_G(g,h)=\max\{d_G(g,h),d_G(g^{-1},h^{-1})\}$, then we have
$$D_G(\Phi(w),\Phi(v))\le\Delta_\Gamma(w,v).$$
It implies that, for any $\Delta_\Gamma$-Cauchy sequence $(w_n)_{n\in\Bbb N}$, $(\Phi(w_n))_{n\in\Bbb N}$ is also a $D_G$-Cauchy sequence.
It is well known that $D_G$ is a compatible complete metric on $G$ (see, e.g. \cite{BK}, Corollary 1.2.2).
Thus $\Phi$ can be uniquely extended to a group homomorphism from ${\overline F}_\Gamma({\cal N}_\omega)$ to $G$.
It is routine to see that, this extend homomorphism, still denoted by $\Phi$, is continuous.

In the end , we show that $\Phi$ is surjection. By (b), for each $g\in G$, there exists an $x\in{\cal N}\subseteq{\overline F}_\Gamma({\cal N}_\omega)$ such that
$$g\in\varphi(\pi_n(x))V_n\cap V_n\varphi(\pi_n(x)).$$
Thus $\varphi(\pi_n(x))^{-1}g\in V_n$, $g\varphi(\pi_n(x))^{-1}\in V_n$, so
$$d_G(g,\varphi(\pi_n(x)))\le 2^{-(n+K+3)},\quad d_G(g^{-1},\varphi(\pi_n(x))^{-1})\le 2^{-(n+K+3)}.$$
It follows that $\Phi(x)=\lim_{n\to\infty}\varphi(\pi_n(x))=g$. Thus $\Phi$ is a surjectively homomorphism from ${\overline F}_\Gamma({\cal N}_\omega)$ to $G$.
This finishes the proof.
\hfill$\Box$

Now we are ready to give two examples of scales $\Gamma$ on $\overline{\cal N}_\omega$ such that $\overline{F}_\Gamma({\cal N}_\omega)$ is surjectively
universal as follows:

\begin{example}
Fix a bijection $\xi:{\cal N}_\omega\to\Bbb N$ with $\xi(\hat 0)=0$. For each $x\in{\cal N}_\omega$ we define
$$\Gamma_1(x^{-1},r)=\Gamma_1(x,r)=(2^6\xi(x)+1)r.$$
\end{example}

Let $K=0$. For any $x\in{\cal N}_m\setminus{\cal N}_{m-1}$, put $y=\hat 0$. We can check that $x\in E_m^k$ implies $\xi(x)<2^{k-(m+5)}$, so $E_m^k$ is finite.
From Theorem~\ref{universal}, we see that $\overline{F}_{\Gamma_1}({\cal N}_\omega)$ is a surjectively universal Polish group.

\begin{example}
Put $\Gamma_2(\hat 0^{-1},r)=\Gamma_2(\hat 0,r)=r$. Fix a bijection $\zeta:{\cal N}_\omega\setminus\{\hat 0\}\to\Bbb N$
such that $x\in{\cal N}_{\zeta(x)+1}$ for each $x$. We define
$$\Gamma_2(x^{-1},r)=\Gamma_2(x,r)=\left\{\begin{array}{ll}8r, & r\le 2^{-(\zeta(x)+6)},\cr \max\{1/8,r\}, & r> 2^{-(\zeta(x)+6)}.\end{array}\right.$$
\end{example}

Let $K=0$. For any $x\in{\cal N}_m\setminus{\cal N}_{m-1}$, put $y=\hat 0$. We can check that $x\in E_m^k$ implies $\zeta(x)<k-6$, so $E_m^k$ is finite.

We conclude this section by:

\begin{theorem}
There exist surjectively universal Polish groups.
\end{theorem}

\section{Adequate scales and good scales}

Let $(X,d)$ be a metric space. Given a scale $\Gamma$ on $\overline{X}$ and two elements $w,v\in F(X)$,
since $N_\Gamma(w^{-1}\cdot v)$ is defined as an infimum, it is almost impossible to compute the exact value of $\delta_\Gamma(w,v)$ directly.
We say a scale on $\overline{X}$ is {\it trivial}, denoted by $\Gamma_{\rm Gr}$, if for all $x\in\overline{X}$ and $r\ge 0$,
$$\Gamma_{\rm Gr}(x,r)=r.$$
We can see that the metric induced by $\Gamma_{\rm Gr}$ is exactly the Graev metric on $F(X)$
(see \cite{graev} or \cite{DG1}, Definition 3.1). It was shown in Theorem 3.6 of \cite{DG1} that, in computing the Graev metric $\delta(w,e)$ for any $w\in F(X)$,
no trivial extension of $w$ is needed. In this section, we try to generalize this method to some widely applicable type of scales.

Now we focus on such a problem: for what type of scales $\Gamma$, whenever we obtain a word $v$,
by eliminating one occurrence of $e$ or $z^{-1}z$ from a given word $w\in W(X)$, for any match $\theta$, can we find another match $\mu$ such that
$$N_\Gamma^\theta(w)\ge N_\Gamma^\mu(v)?$$
If this requirement holds for a given scale $\Gamma$, we can eliminate each occurrence of $e$ or $z^{-1}z$ from $w$, step by step,
to obtain a sequence $w=w_0,w_1,\cdots,w_{n-1},w_n=w'$.
At the same time, we obtain a sequence of matches $\theta_0,\cdots,\theta_n$ from any given match $\theta$ on $\{0,\cdots,{\rm lh}(w)-1\}$ such that
\begin{enumerate}
\item[(i)] $\theta_0=\theta$; for each $i\le n$, $\theta_i$ is a match on $\{0,\cdots,{\rm lh}(w_i)-1\}$; and
\item[(ii)] $N_\Gamma^{\theta_i}(w_i)\ge N_\Gamma^{\theta_{i+1}}(w_{i+1})$ for $i<n$.
\end{enumerate}
Following this method, we can deduce a formula for $u\in F(X)$ as follows:
$$N_\Gamma(u)=\inf\{N_\Gamma^\theta(w):w'=u,\theta\mbox{ is a match}\}=\min\{N_\Gamma^\vartheta(u):\vartheta\mbox{ is a match}\}.$$
Since there are only finitely many matches on $\{0,\cdots,{\rm lh}(u)-1\}$, we can compute $N_\Gamma(u)$ directly.

Given a match $\theta$ on $\{0,\cdots,m\}$, for $0<k<l<m$, it is clear that $\theta\upharpoonright\{k,\cdots,l\}$ is still a match iff
$k\le\theta(i)\le l$ whenever $k\le i\le l$. Let $w=w_0w_1w_2\in W(X)$ with ${\rm lh}(w)=m+1$, and let the domain of $w_1$ be $\{k,\cdots,l\}$.
For future notational simplicity, we make convention that
$$N_\Gamma^\theta(w_1)=N_\Gamma^{\theta\upharpoonright\{k,\cdots,l\}}(w_1).$$

Let $w,v\in W(X)$ with ${\rm lh}(w)=m+1$, ${\rm lh}(v)=n+1$, and let $\theta$ and $\mu$ be two matches on $\{0,\cdots,m\}$ and $\{0,\cdots,n\}$, respectively.
If $w=w_0v_1w_2$, $v=v_0v_1v_2$ with ${\rm lh}(w_0)=k_1$, ${\rm lh}(v_0)=k_2$, ${\rm lh}(v_1)=l$, we say $\theta$ and $\mu$ are {\it coincide} for $v_1$ if
\begin{enumerate}
\item[(i)] both $\theta\upharpoonright\{k_1,\cdots,k_1+l-1\}$ and $\mu\upharpoonright\{k_2,\cdots,k_2+l-1\}$ are still matches; and
\item[(ii)] $\theta(i)=\mu((k_2-k_1)+i)-(k_2-k_1)$ for each $i\in\{k_1,\cdots,k_1+l-1\}$.
\end{enumerate}
Similarly, if $w=v_0w_1v_2$, $v=v_0v_1v_2$ with ${\rm lh}(v_0)=k$ and ${\rm lh}(v_2)=l$.
Denote by $\sigma$ the strictly increasing bijection from $\{0,\cdots,k-1\}\cup\{m-l+1,\cdots,m\}$ to $\{0,\cdots,k-1\}\cup\{n-l+1,\cdots,n\}$.
We say $\theta$ and $\mu$ are {\it coincide} for $(v_0,v_2)$ if
\begin{enumerate}
\item[(i)] both $\theta\upharpoonright\{k,\cdots,m-l\}$ and $\mu\upharpoonright\{k,\cdots,n-l\}$ are still matches; and
\item[(ii)] $\theta(i)=\sigma^{-1}\mu\sigma(i)$ for each $i\in\{0,\cdots,k-1\}\cup\{m-l+1,\cdots,m\}$.
\end{enumerate}

A straightforward observation gives the following fact about pre-norm.

\begin{proposition}\label{segment}
Suppose $\theta$ and $\mu$ are {\it coincide} for $(v_0,v_2)$. Then for any scale $\Gamma$, if $N_\Gamma^\theta(w_1)\ge N_\Gamma^\mu(v_1)$, we have
$$N_\Gamma^\theta(v_0w_1v_2)\ge N_\Gamma^\mu(v_0v_1v_2).$$
\end{proposition}

{\bf Proof.} From Definition~\ref{prenorm}, this will be done by a routine induction.
\hfill$\Box$

When we eliminate one occurrence of $e$ or $z^{-1}z$ with $z\in\overline{X}$, we need to consider all cases of combinations of eliminated alphabets and
effects of a given match on positions of the alphabets. These include eight cases (and their symmetric cases), six trivial and two non-trivial cases, as follows:

Let $\theta$, $\mu$ be matches on $\{0,\cdots,{\rm lh}(w)-1\}$ and $\{0,\cdots,{\rm lh}(v)-1\}$, respectively.

{\sl Trivial case 1.} If $w=v_0ev_1$, $v=v_0v_1$, for brevity, we use figure
$$v_0\;\stackrel{\frown}{e}\;v_1\quad\longrightarrow\quad v_0\;v_1$$
to express that $\theta$ and $\mu$ are coincide for $(v_0,v_1)$, and by letting ${\rm lh}(v_0)=k$, we have $\theta(k)=k$.
By Proposition~\ref{segment}, since $N_\Gamma^\theta(e)=d(e,e)=0$, we have $N_\Gamma^\theta(w)=N_\Gamma^\mu(u)$ for any scale $\Gamma$.

{\sl Trivial case 2.} If $w=v_0ev_1xv_2$, $v=v_0v_1xv_2$, for brevity, we use figure
$$v_0\;\overbrace{e\,v_1x}\;v_2\quad\longrightarrow\quad v_0\;v_1\;\stackrel{\frown}{x}\;v_2$$
to express that $\theta$ and $\mu$ are coincide for both $v_1$ and $(v_0,v_2)$, and by letting ${\rm lh}(v_0)=k$, ${\rm lh}(v_0ev_1)=l$,
we have $\theta(k)=l$, $\mu(l-1)=l-1$. In this case, since
$$\begin{array}{ll}N_\Gamma^\theta(ev_1x)&=\min\{\Gamma(e,N_\Gamma^\theta(v_1)),\Gamma(x,N_\Gamma^\theta(v_1))\}+d(e,x)\cr
&=N_\Gamma^\theta(v_1)+d(e,x)\cr
&=N_\Gamma^\mu(v_1)+d(e,x)=N_\Gamma^\mu(v_1x),\end{array}$$
we have $N_\Gamma^\theta(w)=N_\Gamma^\mu(u)$ for any scale $\Gamma$.

We omit the explanations for figures in following cases.

{\sl Trivial case 3.}
$$v_0\;\stackrel{\frown}{z^{-1}}\stackrel{\frown}{z}\;v_1\quad\longrightarrow\quad v_0\;v_1$$
For any scale $\Gamma$, we have $N_\Gamma^\theta(z^{-1}z)=d(z^{-1},e)+d(z,e)=2d(z,e)\ge 0$.

{\sl Trivial case 4.}
$$v_0\;\overbrace{z^{-1}z}\;v_1\quad\longrightarrow\quad v_0\;v_1$$
For any scale $\Gamma$, we have $N_\Gamma^\theta(z^{-1}z)=0$.

{\sl Trivial case 5.}
$$v_0\;\stackrel{\frown}{z}\overbrace{z^{-1}\,v_1x}\;v_2\quad\longrightarrow\quad v_0\;v_1\;\stackrel{\frown}{x}\;v_2$$
For any scale $\Gamma$, we have
$$\begin{array}{ll}N_\Gamma^\theta(zz^{-1}v_1x)&=d(z,e)+\min\{\Gamma(z,N_\Gamma^\theta(v_1)),\Gamma(x,N_\Gamma^\theta(v_1))\}+d(z,x)\cr
&\ge d(z,e)+N_\Gamma^\theta(v_1)+d(z,x)\cr
&\ge N_\Gamma^\mu(v_1)+d(x,e)=N_\Gamma^\mu(v_1x).\end{array}$$

{\sl Trivial case 6.}
$$v_0\;\overbrace{z^{-1}\stackrel{\frown}{z}v_1x}\;v_2\quad\longrightarrow\quad v_0\;v_1\;\stackrel{\frown}{x}\;v_2$$
For any scale $\Gamma$, we have
$$\begin{array}{ll}N_\Gamma^\theta(z^{-1}zv_1x)&=\min\{\Gamma(z,d(z,e)+N_\Gamma^\theta(v_1)),\Gamma(x,d(z,e)+N_\Gamma^\theta(v_1))\}+d(z,x)\cr
&\ge d(z,e)+N_\Gamma^\theta(v_1)+d(z,x)\cr
&\ge N_\Gamma^\mu(v_1)+d(x,e)=N_\Gamma^\mu(v_1x).\end{array}$$

Unlike the preceding six trivial cases, the following two non-trivial cases is true restriction on the scale under consideration.

{\sl Non-trivial case 1.}
$$v_0\;\overbrace{x^{-1}v_1z}\;\overbrace{z^{-1}v_2y}\;v_3\quad\longrightarrow\quad v_0\;\overbrace{x^{-1}v_1v_2y}\;v_3$$
Denote $r_1=N_\Gamma^\theta(v_1)$, $r_2=N_\Gamma^\theta(v_2)$. For this case, we need the following inequality holds for $\Gamma$.
\begin{enumerate}
\item[(A1)]
$\min\{\Gamma(x,r_1),\Gamma(z,r_1)\}+d(x,z)+\min\{\Gamma(y,r_2),\Gamma(z,r_2)\}+d(y,z)\\
\ge\min\{\Gamma(x,r_1+r_2),\Gamma(y,r_1+r_2)\}+d(x,y)$.
\end{enumerate}

{\sl Non-trivial case 2.}
$$v_0\;\overbrace{x^{-1}v_1\;\underbrace{y\,v_2z}\!\,^{-1}\;z}\;v_3\quad\longrightarrow\quad v_0\;\overbrace{x^{-1}v_1y}\;v_2\;v_3$$
Denote $r_1=N_\Gamma^\theta(v_1)$, $r_2=N_\Gamma^\theta(v_2)$. Note that $d(y^{-1},z^{-1})=d(y,z)$.
For this case, we need the following inequality holds for $\Gamma$.
\begin{enumerate}
\item[(A2)] Let $r=r_1+\min\{\Gamma(y^{-1},r_2),\Gamma(z^{-1},r_2)\}+d(y,z)$, then
$$\min\{\Gamma(x,r),\Gamma(z,r)\}+d(x,z)\ge\min\{\Gamma(x,r_1),\Gamma(y,r_1)\}+d(x,y)+r_2.$$
\end{enumerate}

\begin{definition}{\rm
Let $\Gamma$ be a scale on $\overline{X}$. We say $\Gamma$ is {\it adequate} if, for any $x,y,z\in\overline{X}$ and $r_1,r_2\ge 0$,
conditions (A1) and (A2) hold for $\Gamma$.}
\end{definition}

From the previous arguments, we get the following theorem.

\begin{theorem}\label{compute}
Let $\Gamma$ be an adequate scale on $\overline{X}$. For any $w\in F(X)$, we have
$$N_\Gamma(w)=\min\{N_\Gamma^\theta(w):\theta\mbox{ is a match}\}.$$
\end{theorem}

An obvious but useful corollary of Theorem~\ref{compute} is the following:

\begin{corollary}\label{sub}
Let $\Gamma$ be an adequate scale on $\overline{X}$ and $Y\subseteq X$. Let $\delta_\Gamma^Y$ be the metric induced from
the subspace $(\overline{Y},d\upharpoonright\overline{Y})$ with the scale $\Gamma\upharpoonright\overline{Y}$.
Then $\delta_\Gamma^Y=\delta_\Gamma\upharpoonright\overline{Y}$. Moreover, if $Y$ is dense in $X$, then $\overline{F}_\Gamma(Y)=\overline{F}_\Gamma(X)$.
\end{corollary}

It is obvious that the trivial scale $\Gamma_{\rm Gr}(x,r)\equiv r$ is adequate, so Theorem 3.6 of \cite{DG1} is a corollary of this theorem.

Conditions (A1) and (A2) are so complicated that it is very hard to check them for a given scale. We are going to simplify these conditions.

\begin{lemma}
Condition {\rm (A1)} is equivalent to the following:
\begin{enumerate}
\item[(A1)$'$] $\Gamma(z,r)+d(x,z)+d(y,z)\ge\min\{\Gamma(x,r),\Gamma(y,r)\}+d(x,y)$; and
\item[(A1)$''$] $\min\{\Gamma(x,r_1),\Gamma(y,r_1)\}+\min\{\Gamma(x,r_2),\Gamma(y,r_2)\}\\ \ge\min\{\Gamma(x,r_1+r_2),\Gamma(y,r_1+r_2)\}$.
\end{enumerate}
\end{lemma}

{\bf Proof.} For (A1)$\Rightarrow$(A1)$'$, put $r=r_1$, $r_2=0$.

For (A1)$\Rightarrow$(A1)$''$, putting $y=z$ in (A1), we get
$$\min\{\Gamma(x,r_1),\Gamma(y,r_1)\}+\Gamma(y,r_2)\ge\min\{\Gamma(x,r_1+r_2),\Gamma(y,r_1+r_2)\}.$$
By changing position of $x$ and $y$, we get
$$\min\{\Gamma(x,r_1),\Gamma(y,r_1)\}+\Gamma(x,r_2)\ge\min\{\Gamma(x,r_1+r_2),\Gamma(y,r_1+r_2)\}.$$
Hence (A1)$''$ follows.

For (A1)$'$+(A2)$''$ $\Rightarrow$(A1). Firstly, we have
$$\begin{array}{ll} & \min\{\Gamma(x,r_1+r_2),\Gamma(y,r_1+r_2)\}+d(x,y)\cr
\le & \Gamma(z,r_1+r_2)+d(x,z)+d(y,z)\cr
\le & \Gamma(z,r_1)+\Gamma(z,r_2)+d(x,z)+d(y,z).\end{array}$$
Secondly, we have
$$\begin{array}{ll} & \min\{\Gamma(x,r_1+r_2),\Gamma(y,r_1+r_2)\}+d(x,y)\cr
\le & \min\{\Gamma(x,r_1),\Gamma(y,r_1)\}+\min\{\Gamma(x,r_2),\Gamma(y,r_2)\}+d(x,y)\cr
\le & \min\left\{\begin{array}{l}\Gamma(x,r_1)+\Gamma(y,r_2)+d(x,z)+d(y,z)\cr \Gamma(x,r_1)+\Gamma(z,r_2)+d(x,z)+d(y,z)\cr \Gamma(z,r_1)+d(x,z)+d(y,z)+\Gamma(y,r_2)
\end{array}\right\}.\end{array}$$
From these (A1) follows.
\hfill$\Box$

It is not easy to simplify condition (A2). Condition (A1)$''$ is still too complicated. We turn to find some sufficient conditions for (A1)$''$ and (A2).

\begin{lemma}\label{g23}
Let $\Gamma$ be a scale on $\overline{X}$.
\begin{enumerate}
\item[(i)] If $\Gamma(x,r)/r$ is monotone decreasing with respect to variable $r$, then {\rm (A1)}$''$ holds.
\item[(ii)] If {\rm (A1)}$'$ holds, and for any $x\in\overline{X}$, $r_1,r_2\ge 0$, we have $\Gamma(x,r_1+r_2)\ge\Gamma(x,r_1)+r_2$, then {\rm (A2)} holds.
\end{enumerate}
\end{lemma}

{\bf Proof.} (i) Let $x,y\in\overline{X}$. Define function $D:\Bbb R_+\to\Bbb R_+$ as
$$D(r)=\min\{\Gamma(x,r)/r,\Gamma(y,r)/r\}.$$
Then $D$ is a monotone decreasing function. Thus
$$\begin{array}{ll}&\min\{\Gamma(x,r_1),\Gamma(y,r_1)\}+\min\{\Gamma(x,r_2),\Gamma(y,r_2)\}\cr
=&r_1D(r_1)+r_2D(r_2)\cr
\ge &(r_1+r_2)D(r_1+r_2)\cr
=&\min\{\Gamma(x,r_1+r_2),\Gamma(y,r_1+r_2)\}.\end{array}$$

(ii) Note that
$$r=r_1+\min\{\Gamma(y^{-1},r_2),\Gamma(z^{-1},r_2)\}+d(y,z)\ge r_1+r_2+d(y,z).$$
Therefore, we have
$$\begin{array}{ll}&\min\{\Gamma(x,r),\Gamma(z,r)\}+d(x,z)\cr
\ge &\min\{\Gamma(x,r_1+r_2+d(y,z)),\Gamma(z,r_1+r_2+d(y,z))\}+d(x,z)\cr
\ge &\min\{\Gamma(x,r_1),\Gamma(z,r_1)\}+r_2+d(y,z)+d(x,z).\end{array}$$
On the one hand,
$$\Gamma(x,r_1)+r_2+d(y,z)+d(x,z)\ge\min\{\Gamma(x,r_1),\Gamma(y,r_1)\}+d(x,y)+r_2$$
is trivial. On the other hand,
$$\Gamma(z,r_1)+r_2+d(y,z)+d(x,z)\ge\min\{\Gamma(x,r_1),\Gamma(y,r_1)\}+d(x,y)+r_2$$
follows from (A1)$'$. Hence (A2) holds for $\Gamma$.
\hfill$\Box$

Most of the time, we are concerned about scales on such metric spaces $(X,d)$ whose metric $d$ is actually an {\it ultrametric}, i.e.
$$d(x,y)\le\max\{d(x,z),d(y,z)\}.$$
We still extend the ultrametric $d$ to an ultrametric on $\overline{X}$.

\begin{definition}{\rm
Let $(\overline{X},d)$ be an ultrametric space. We say a scale $\Gamma$ is {\it good} if, for any $x,y\in\overline{X}$ and $r,r_1,r_2\in\Bbb R_+$, we have
\begin{enumerate}
\item[(G1)] $\Gamma(y,r)+d(x,y)\ge\Gamma(x,r)$;
\item[(G2)] $\Gamma(x,r)/r$ is monotone decreasing with respect to variable $r$;
\item[(G3)] $\Gamma(x,r_1+r_2)\ge\Gamma(x,r_1)+r_2$.
\end{enumerate}}
\end{definition}

\begin{theorem}
Let $(\overline{X},d)$ be an ultrametric space. If $\Gamma$ is a good scale on $\overline{X}$, then $\Gamma$ is adequate.
\end{theorem}

{\bf Proof.} Firstly, we show that (G1) implies (A1)$'$. For $x,y,z\in\overline{X}$, assume that $d(x,z)\le d(y,z)$.
Since $d$ is an ultrametric, we have $d(x,y)\le d(y,z)$, so
$$\Gamma(z,r)+d(x,z)+d(y,z)\ge\Gamma(x,r)+d(y,z)\ge\min\{\Gamma(x,r),\Gamma(y,r)\}+d(x,y).$$
Then the theorem follows from Lemma~\ref{g23}
\hfill$\Box$

In the end of this section, we give an example of good scale $\Gamma$ on ${\cal N}_\omega$ such that $\overline{F}_\Gamma({\cal N}_\omega)$
is surjectively universal.

\begin{example}
Put $\Gamma_0(\hat 0,r)=r$. Let $m\ge 1$. Suppose we have defined $\Gamma_0(x,r)$ for $x\in{\cal N}_{m-1}$.
We extend the definition onto $x\in{\cal N}_m\setminus{\cal N}_{m-1}$ as
$$\Gamma_0(x^{-1},r)=\Gamma_0(x,r)=\min\left
\{\begin{array}{l}2^5(1+2^{x(0)}+\cdots+2^{x(0)+\cdots+x(m-1)})r\cr \Gamma_0(\pi_{m-1}(x),r)+2^{-m}\end{array}\right\}.$$
\end{example}

Let $K=1$. For any $x\in{\cal N}_m\setminus{\cal N}_{m-1}$, put $y=\pi_{m-1}(x)$. Then $2^{-K}d(x,y)=2^{-m}$.
We can prove that $x\in E_m^k$ implies that, for some $2^{-k}<r\le 2^{-(m+6)}$,
$$\begin{array}{ll}\Gamma_0(x,r)&=2^5(1+2^{x(0)}+\cdots+2^{x_0+\cdots+x(m-1)})r\cr
&\le \Gamma_0(\pi_{m-1}(x),r)+2^{-m}\cr
&\le 2^5(1+2^{x(0)}+\cdots+2^{x_0+\cdots+x(m-2)})r+2^{-m}.\end{array}$$
It follows that $r\le 2^{-(m+5+x(0)+\cdots+x(m-1))}$. Comparing with $2^{-k}<r$, we have
$$x(0)+\cdots+x(m-1)<k-(m+5),$$
so $E_m^k$ is finite. From Theorem~\ref{universal}, $\overline{F}_{\Gamma_0}({\cal N}_\omega)$ is surjectively universal.

To see that $\Gamma_0$ is good, let $x,y\in{\cal N}_\omega$ with $d(x,y)=2^{-n}$, then $\pi_n(x)=\pi_n(y)$. Note that
$$\Gamma_0(\pi_n(x),r)\le\Gamma_0(x,r)\le\Gamma_0(\pi_n(x),r)+\sum_{m>n}2^{-m}<\Gamma_0(\pi_n(x),r)+2^{-n},$$
$$\Gamma_0(\pi_n(x),r)=\Gamma_0(\pi_n(y),r)\le\Gamma_0(y,r)<\Gamma_0(\pi_n(x),r)+2^{-n}.$$
Thus $\Gamma_0(y,r)+2^n\ge\Gamma_0(x,r)$. Then (G1) holds. (G2) and (G3) follow from a routine induction.

Moreover, $\Gamma_0$ can be extend to be a regular scale. Recall that a scale on $\overline{\cal N}$ is {\it regular} if for all $x\in\overline{\cal N}$,
$r\in\Bbb R_+$ and $n\in{\Bbb N}$, we have $\Gamma(x,r)\ge\Gamma(\pi_n(x),r)$ (see \cite{DG},
Definition 4.1). By letting $\Gamma_0(x,r)=\lim_{n\to\infty}\Gamma_0(\pi_n(x),r)$,
we extend $\Gamma_0$ onto $\overline{\cal N}$. The regularity of extended scale is easy to check. Then from Theorem 4.5 of \cite{DG} and Corollary~\ref{sub},
$\overline{F}_{\Gamma_0}({\cal N}_\omega)$ is a ${\bf\Pi}^0_3$ subgroup of the inverse limit $\displaystyle\lim_{\stackrel{\longleftarrow}{n}}F({\cal N}_n)$.
Furthermore, it is isomorphic to a ${\bf\Pi}^0_3$ subgroup of $S_\infty$.

\section{Adequate scales and CLI groups}

In this section, we prove that any non-trivial adequate scale can not induce a CLI group.
Recall that a CLI group is a Polish group admitting a compatible complete left-invariant metric. It is well known that, if $G$ is a CLI group,
then any compatible left-invariant metric on $G$ must be complete (see \cite{gaobook}, Proposition 2.2.6).

\begin{lemma}\label{value}
Let $m,n\in\Bbb N$. If $0<m<2^n$, then there is a $u\in F({\cal N}_\omega)$ such that $N_\Gamma(u)=m/2^n$ for any scale $\Gamma$ on ${\cal N}_\omega$.
\end{lemma}

{\bf Proof.} We can find $0<n_1<\cdots<n_k\le n$ such that
$$m/2^n=2^{-n_1}+\cdots+2^{-n_k}.$$
Select $x_1,y_1,\cdots,x_k,y_k$ from ${\cal N}_\omega$ such that, for $1\le i\le k$,
$$d(x_i,y_i)=2^{-n_i},\quad x_i(0)=y_i(0)=i.$$
Define $u=x_1^{-1}y_1\cdots x_k^{-1}y_k$. Let $\theta_0$ be the match on $\{0,\cdots,2k-1\}$ with $\theta_0(2i)=2i+1$ for $0\le i<k$.
Then for any scale $\Gamma$ we have
$$N_\Gamma(u)\le N_\Gamma^{\theta_0}(u)=d(x_1,y_1)+\cdots+d(x_k,y_k)=m/2^n.$$
Let $\theta\ne\theta_0$ be any match on $\{0,\cdots,2k-1\}$. There exists an $i<k$ such that $\theta(2i)\ne 2i+1$.
Note that $d(x_i,e)=d(x_i,x_j)=d(x_i,x_j^{-1})=1$ for any $j\ne i$, we have $N_{\Gamma_{\rm Gr}}^\theta(u)\ge 1$. It follows from $m/2^n<1$ that
$$N_\Gamma(u)\ge N_{\Gamma_{\rm Gr}}(u)=\min\{ N_{\Gamma_{\rm Gr}}^\theta(u):\theta\mbox{ is a match}\}=m/2^n.$$
Thus $N_\Gamma(u)=m/2^n$.
\hfill$\Box$

\begin{theorem}
Let $\Gamma$ be a non-trivial adequate scale on ${\cal N}_\omega$. Then $\overline{F}_\Gamma({\cal N}_\omega)$ is not CLI group.
\end{theorem}

{\bf Proof.} Since $\Gamma$ is non-trivial, there is an $x_0\in\overline{\cal N}_\omega$ such that $\Gamma(x_0,r)\not\equiv r$.
By the preceding remarks, we assume for contradiction that $\delta_\Gamma$ is a complete metric on $\overline{F}_\Gamma({\cal N}_\omega)$.

Denote $f(r)=\Gamma(x_0,r)$. Since $\Gamma$ is adequate, for $r_1,r_2\ge 0$, (A1)$''$ gives
$$\Gamma(x_0,r_1)+\Gamma(x_0,r_2)\ge\Gamma(x_0,r_1+r_2).\eqno(*)$$
Note that $\Gamma(x_0,r)$ is monotone increasing and $\lim_{r\to 0}\Gamma(x_0,r)=0$, we see that $f(r)$ is continuous.

Firstly, we claim that there is $a\in(0,1)$ such that $f(r)>r$ for all $0<r\le a$. Assume for contradiction that there exists a sequence of real numbers
$a_1>a_2>\cdots>a_k>\cdots>0$ such that $\lim_{k\to\infty}a_k=0$ and $f(a_k)=a_k$ for each $k$.
Then inequality $(*)$, together with $\Gamma(x_0,r)\ge r$ and continuity of $f(r)$, gives $\Gamma(x_0,r)\equiv r$. A contradiction!

Secondly, we claim that, for any $0<r<a$, $f^k(r)\ge a$ for large enough $k$, in which $f^k=\overbrace{f\circ\cdots\circ f}^k$.
Otherwise, $(f^k(r))_{k\in\Bbb N}$ forms a strictly increasing sequence below $a$, thus it converge to some $b\le a$. By continuity, we have $f(b)=b$.
A contradiction!

Without loss of generality, assume that $x_0(0)=0$ and $a\le 1/2$.
Find a sequence of natural numbers $0=k_0<k_1<\cdots<k_j<\cdots$ such that $f^{k_j}(2^{-(j+1)})\ge a$ for $j\ge 0$.
For each $j$, by continuity of $f$, there are $m_j,n_j\in\Bbb N$ with $0<m_j<2^{n_j}$ such that $2^{-(j+1)}\le f^{k_{j+1}-k_j}(m_j/2^{n_j})<2^{-j}$.
From Lemma~\ref{value}, we can find $u_{2j},u_{2j+1}\in F(X)$ such that $N_\Gamma(u_{2j})=m_j/2^{n_j}$ and $N_\Gamma(u_{2j+1})=2^{-(j+2)}$.
Furthermore, we can assume that every alphabet $x$ appearing in each $u_l$ satisfies that $x(0)\ne 0$, thus $d(x_0,x)=1$. Let
$$v_l=\left\{\begin{array}{ll}(x_0^{-1})^{k_{j+1}-k_j}u_{2j}x_0^{k_{j+1}-k_j}, & l=2j,\cr u_{2j+1}, & l=2j+1,\end{array}\right.$$
and $w_l=v_0v_1\cdots v_l$ for $l\in\Bbb N$. We can see that all $w_l$ are irreducible words.

On the one hand, it is clear that $N_\Gamma(v_{2j})\le f^{k_j-k_{j-1}}(m_j/2^{n_j})<2^{-j}$ and $N_\Gamma(v_{2j+1})=2^{-(j+2)}$. It follows that
$$\sum_{l\ge 1}\delta_\Gamma(w_{l-1},w_l)=\sum_{l\ge 1}N_\Gamma(w_{l-1}^{-1}\cdot w_l)=\sum_{l\ge 1}N_\Gamma(v_l)<+\infty.$$
So $(w_l)_{l\in\Bbb N}$ is a $\delta_\Gamma$-Cauchy sequence. Since $\delta_\Gamma$ is complete metric,
there is $\sigma\in\overline{F}_\Gamma({\cal N}_\omega)$ such that $\lim_{l\to\infty}w_l=\sigma$. Thus $\lim_{l\to\infty}w_l^{-1}=\sigma^{-1}$.

On the other hand, given an $l\ge 1$, the following is an irreducible word:
$$w_l\cdot w_{l-1}^{-1}=v_0v_1\cdots v_{l-1}v_lv_{l-1}^{-1}\cdots v_1^{-1}v_0^{-1}.$$
We rewrite $w_l\cdot w_{l-1}^{-1}=z_0\cdots z_n$ with $z_0,\cdots,z_n\in\overline{\cal N}_\omega$. Denote
$$I_0=\{i\le n:z_i=x_0\},\quad I_1=\{i\le n:z_i=x_0^{-1}\},\quad I=I_0\cup I_1.$$

Let $\theta$ be a match on $\{0,\cdots,n\}$. If there is an $i\in I_0$ such that $\theta(i)\notin I_1$,
then it is easy to see that $N_\Gamma^\theta(w_l\cdot w_{l-1}^{-1})\ge 1$. Now assume that $i\in I_0$ iff $\theta(i)\in I_1$.
Let $i_0$ be the least $i\in I$ such that $\theta(i)<i$. Note that $\theta\upharpoonright\{\theta(i_0)+1,\cdots,i_0-1\}$ is still a match.
Comparing with the definition of $i_0$, we see that $I\cap\{\theta(i_0)+1,\cdots,i_0-1\}=\emptyset$.

{\sl Case 1.} If $i_0\in I_0$, then $z_{i_0}=x_0$, $z_{\theta(i_0)}=x_0^{-1}$. This implies that there is a $j_0$ such that
$u_{2j_0}=z_{\theta(i_0)+1}\cdots z_{i_0-1}$. There are $k_{j_0+1}$ many $x_0^{-1}$'s at the left side of $u_{2j_0}$, i.e.
the set $\{i\in I_1:i<i_0\}$ has exact $k_{j_0+1}$ many elements. Denote $k=k_{j_0+1}$. We enumerate $\{i\in I_1:i<i_0\}$ as
$$i_1<\cdots<i_{k-1}<i_k=\theta(i_0)<i_0.$$
From the definition of $i_0$, we see that $\theta(i_1),\cdots,\theta(i_{k-1})\ge i_0$. The definition of match gives
$$\theta(i_1)>\cdots>\theta(i_{k-1})>\theta(i_k)=i_0.$$
By repeatedly applying Proposition~\ref{segment}, we get
$$\begin{array}{ll}N_\Gamma^\theta(w_l\cdot w_{l-1}^{-1})&\ge f^k(N_\Gamma^\theta(u_{2j_0}))\ge f^{k_{j_0}}(f^{k_{j_0+1}-k_{j_0}}(N_\Gamma(u_{2j_0})))\cr
&=f^{k_{j_0}}(f^{k_{j_0+1}-k_{j_0}}(m_{j_0}/2^{n_{j_0}}))\cr
&\ge f^{k_{j_0}}(2^{-(j_0+1)})\ge a.\end{array}$$

{\sl Case 2.} If $i_0\in I_1$, then $z_{i_0}=x_0^{-1}$, $z_{\theta(i_0)}=x_0$. This implies that there is a $j_0\ge 1$ such that
$u_{2j_0-1}=z_{\theta(i_0)+1}\cdots z_{i_0-1}$. There are $k_{j_0}$ many $x_0^{-1}$'s at the left side of $u_{2j_0-1}$.
By similar arguments as in case 1, we can prove that $N_\Gamma^\theta(w_l\cdot w_{l-1}^{-1})\ge a$.

Since $\Gamma$ is adequate, for every $l\ge 1$, we have
$$\delta_\Gamma(w_l^{-1},w_{l-1}^{-1})=N_\Gamma(w_l\cdot w_{l-1}^{-1})=\min\{N_\Gamma^\theta(w_l\cdot w_{l-1}^{-1}):\theta\mbox{ is a match}\}\ge a.$$
This contradicts with $\lim_{l\to\infty}w_l^{-1}=\sigma^{-1}$.
\hfill$\Box$



\subsection*{Acknowledgments}

I am grateful to Kechris for a comment concerning the introduction.
The paper was written during the author's visit to the Research Training Group in Logic and Dynamics at the University of North Texas.
I would like to thank UNT and the RTG for the hospitality and travel support.
I also thank Steve Jackson and Mingzhi Xuan for the conversations in seminars during the visit.
Special thanks are due to Su Gao for inviting me to visit UNT and for useful suggestions on the paper.

\end{document}